\documentclass[11pt]{article} 
\usepackage[spanish,english]{babel}
\usepackage{amsmath}
\usepackage{amsfonts}
\usepackage{latexsym}
\usepackage{epsfig}
\usepackage{graphicx}
\usepackage{url}

\newtheorem{theorem}{Theorem}[section]

\newtheorem{proposition}[theorem]{Proposition}
\newtheorem{corollary}[theorem]{Corollary}
\newcommand\qedsymbol{\hbox{\rlap{$\sqcap$}$\sqcup$}}
\newcommand\qed{\relax\ifmmode\else\hbox{}\unskip\fi\hfill\qedsymbol}

\begin{document}
\selectlanguage{english}

\title{$0/1$-Polytopes related to Latin squares autotopisms.\footnote{Official printed version available in Proceedings of VI Jornadas de Matem\'atica Discreta y Algor\'\i tmica (2008), pp. 311-319. ISBN: 978-84-8409-263-6.}}

\author{R. M. Falc\'on\\
Department of Applied Mathematic I. \\
Technical Architecture School. University of Seville. \\
Avda. Reina Mercedes, 4A - 41012, Seville (Spain). \\
\url{rafalgan@us.es}
}

\maketitle

\begin{abstract}
The set $LS(n)$ of Latin squares of order $n$ can be represented in $\mathbb{R}^{n^3}$ as a $(n-1)^3$-dimensional $0/1$-polytope. Given an autotopism $\Theta=(\alpha,\beta,\gamma)\in\mathfrak{A}_n$, we study in this paper the $0/1$-polytope related to the subset of $LS(n)$ having $\Theta$ in their autotopism group. Specifically, we prove that this polyhedral structure is generated by a polytope in $\mathbb{R}^{((\mathbf{n}_{\alpha}-\mathbf{l}_{\alpha}^1)\cdot n^2 + \mathbf{l}_{\alpha}^1\cdot \mathbf{n}_{\beta}\cdot n)-(\mathbf{l}_{\alpha}^1\cdot \mathbf{l}_{\beta}^1\cdot (n -\mathbf{l}_{\gamma}^1) + \mathbf{l}_{\alpha}^1\cdot \mathbf{l}_{\gamma}^1\cdot (\mathbf{n}_{\beta} -\mathbf{l}_{\beta}^1) +
\mathbf{l}_{\beta}^1\cdot \mathbf{l}_{\gamma}^1\cdot (\mathbf{n}_{\alpha} -\mathbf{l}_{\alpha}^1))}$, where $\mathbf{n}_{\alpha}$ and $\mathbf{n}_{\beta}$ are the number of cycles of $\alpha$ and $\beta$, respectively, and $\mathbf{l}_{\delta}^1$ is the number of fixed points of $\delta$, for all $\delta\in \{\alpha,\beta,\gamma\}$. Moreover, we study the dimension of these two polytopes for Latin squares of order up to $9$.
\end{abstract}

{\bf Key words:} $0/1$-polytope, Latin Square, Autotopism group.

\section{Introduction}

A {\em $0/1$-polytope} \cite{Ziegler00} in $\mathbb{R}^d$ is the convex hull $\mathcal{P}$ of a finite set of points with $0/1$-coordinates. Equivalently, it is a polytope with all its vertices in the vertex set of the unit cube $C_d=[0,1]^d$. Thus, if we consider these vertices as the column vectors of a matrix $V\in \{0,1\}^{d\times n}$, it is verified that $\mathcal{P}=\mathcal{P}(V)=conv(V)=\{V\cdot (x_1,x_2,...,x_n)^t \mid x_i\geq 0, \forall i\in [n] \text{ and } \sum_{i\in [n]}x_i=1\},$ where $[n]$ will denote from now on the set $\{1,2,...,n\}$. The {\em dimension} of $\mathcal{P}$ is the maximum number of affinely independent points in $\mathcal{P}$ minus $1$. Permuting coordinates and {\em switching} (replacing $x_i$ by $1-x_i$) coordinates transform $0/1$-polytopes into $0/1$-polytopes. Two $0/1$-polytopes are said to be $0/1$-equivalent if there exists a sequence of the two previous operations transforming one of them into the other one. In combinatorial optimization there are several examples of $0/1$-polytopes like the salesman polytope \cite{Grotschel88}, the cut polytope \cite{Deza97} or the Latin square polytope \cite{Euler86}. In this paper, we are interested in the last one, which appears in the {\em 3-dimensional planar assignment problem} ($3PAP_n$):{\small
$$\min \sum_{i\in I, j\in J, k\in K} w_{ijk}\cdot x_{ijk},\ s.t. \begin{cases}\begin{array}{lcr}
         \sum_{i\in I}x_{ijk}=1, \forall j\in J, k\in K. & \ & (1.1)\\
         \sum_{j\in J}x_{ijk}=1, \forall i\in I, k\in K. & \ & (1.2)\\
         \sum_{k\in K}x_{ijk}=1, \forall i\in I, j\in J. & \ & (1.3)\\
         x_{ijk}\in \{0,1\}, \forall i\in I, j\in J, k \in K. & \ & (1.4)\end{array}
       \end{cases} \hspace{0.05cm} (1)
$$}
where $w_{ijk}$ are real weights and $I,J,K$ are three disjoint $n$-sets.

Euler et al. \cite{Euler86} observed that there exists a 1-1 correspondence between the set $LS(n)$ of Latin squares of order $n$ and the set $FS(n)$ of feasible solutions of the $3PAP_n$. Specifically, a {\em Latin square} $L$ of order $n$ is an $n \times n$ array with elements chosen from a set of $n$ distinct symbols such that each symbol occurs precisely once in each row and each column. From now on, we will assume $[n]$ as this set of symbols. Given $L=\left(l_{i,j}\right)\in LS(n)$, the {\em orthogonal array representation of} $L$ is the set of $n^2$ triples $\{(i,j,l_{i,j})\,\mid \, i,j\in [n]\}$. So, by taking $I=J=K=[n]$ and by considering the lexicographical order in $I\times J\times K$, it can be defined the 1-1 correspondence $\Phi:LS(n)\rightarrow FS(n)\subseteq \mathbb{R}^{n^3}$, such that, given $L=(l_{i,j})\in LS(n)$, it is $\Phi(L)=(x_{111},x_{112},...,x_{1nn},x_{211},...,x_{nnn})$, where $x_{ijk}=\begin{cases} 1, \text{ if } l_{i,j}=k,\\ 0, \text{ otherwise}.\end{cases}$. Moreover, if $A$ is the constraint matrix of the system of equations $(1)$, it is defined the {\em Latin square polytope}, $\mathcal{P}_{LS(n)}=conv\{FS(n)\}=conv\{{\bf x}\in \{0,1\}^{n^3} \mid A\cdot {\bf x}={\bf e}\}$, where ${\bf e}=(1,...,1)^t$ with $3\cdot n^2$ entries. Thus, every point of $\mathcal{P}_{LS}\cap C_{n^3}$ is a Latin square of order $n$ and vice versa. By obtaining the minimal equation system for $P_{LS}$, Euler et al. proved that this polytope is $(n-1)^3$-dimensional and they gave some general results about its facial structure.

In this paper, we are interested in obtaining a similar construction than the above one, in the case of adding some extra conditions to the $3PAP_n$. Specifically, we want to study those $0/1$-polytopes related to Latin squares having some symmetrical restrictions. To expose the problem, some previous considerations are needed: The permutation group on $[n]$ is denoted by $S_n$. Every permutation $\delta\in S_n$ can be uniquely written as a
composition of pairwise disjoint cycles, $\delta=C^{\delta}_1\circ C^{\delta}_2\circ ... \circ
C^{\delta}_{\mathbf{n}_{\delta}},$ where for all $i\in [\mathbf{n}_{\delta}]$, one has
$C^{\delta}_i=\left(c_{i,1}^{\delta}\ c_{i,2}^{\delta}\ ...\
c_{i,\ \lambda_i^{\delta}}^{\delta}\right)$, with
$c_{i,1}^{\delta}=\min_j \{c_{i,j}^{\delta}\}$. The {\em cycle
structure of} $\delta$ is the sequence
$\mathbf{l}_{\delta}=(\mathbf{l}_1^{\delta},\mathbf{l}_2^{\delta},...,\mathbf{l}_n^{\delta})$,
where $\mathbf{l}_i^{\delta}$ is the number of cycles of length
$i$ in $\delta$, for all $i\in [n]$.  Thus, $\mathbf{l}_1^{\delta}$ is the cardinal of the set of {\em fixed points} of $\delta$, $Fix(\delta)=\{i\in [n] \,\mid \, \delta(i)=i\}$. An {\em isotopism} of a Latin square $L=\left(l_{i,j}\right)\in LS(n)$ is a triple
$\Theta=(\alpha,\beta,\gamma)\in \mathcal{I}_n=S_n\times S_n\times
S_n$. In this way, $\alpha,\beta$ and $\gamma$ are permutations of
rows, columns and symbols of $L$, respectively. The resulting
square $L^{\Theta}=\{(\alpha(i),\beta(j),\gamma\left(l_{i,j}\right))\,\mid
\, i,j\in [n]\}$ is also a Latin square. The {\em cycle structure} of $\Theta$ is the triple
$(\mathbf{l}_{\alpha},\mathbf{l}_{\beta},\mathbf{l}_{\gamma})$.

An isotopism which maps $L$ to itself is an {\em autotopism}. The stabilizer subgroup of $L$ in $\mathcal{I}_n$
is its {\em autotopism group}, $\mathfrak{A}(L)=\{\Theta\in
\mathcal{I}_n\,\mid \,L^{\Theta}=L\}$. The set of all autotopisms of Latin squares of order $n$ is denoted by $\mathfrak{A}_n$. Given $\Theta\in\mathfrak{A}_n$, the set of all Latin squares $L$ such
that $\Theta\in \mathfrak{A}(L)$ is denoted by $LS(\Theta)$ and
the cardinality of $LS(\Theta)$ is denoted by $\Delta(\Theta)$.
Specifically, if $\Theta_1$ and $\Theta_2$ are two autotopisms
with the same cycle structure, then
$\Delta(\Theta_1)=\Delta(\Theta_2)$. The
possible cycle structures of the set of non-trivial autotopisms of
Latin squares of order up to $11$ were obtained in
\cite{FalconAC}.

Gr\"obner bases were used in \cite{FalconMartinJSC07} to describe
an algorithm that allows one to obtain the number $\Delta(\Theta)$
in a computational way. This algorithm was implemented in {\sc
Singular} \cite{Greuel05} to get the number of Latin squares of
order up to $7$ related to any autotopism of a given cycle
structure. Specifically, the authors
followed the ideas implemented by Bayer \cite{Bayer82} to solve the problem of an
$n$-colouring a graph, since every Latin square of order $n$ is
equivalent to an $n$-coloured bipartite graph $K_{n,n}$. More recently, Falc\'on and Mart\'\i n-Morales \cite{FalconMartinEACA08} have studied the case $n>7$ by implementing in a new algorithm the 1-1 correspondence between the $3PAP_n$ and the set $LS(n)$. As an immediate consequence, the set of vertices of $C_{n^3}$ related to $LS(\Theta)$ can be obtained.

In Section 2, given $\Theta\in \mathfrak{A}_n$, we study the set of constraints which can be added to the $3PAP_n$ to get a set of feasible solutions equivalent to the set $LS(\Theta)$. In Section 3, we define the $0/1$-polytope in $\mathbb{R}^{n^3}$ related to $LS(\Theta)$. Moreover, we prove the existence of a $0/1$-subpolytope of the previous one which can generate it. We see that these two polytopes do not depend on the autotopism $\Theta$ but on the cycle structure of the autotopism. Finally, we study the dimensions of these polytopes and we give a classification for polytopes related to autotopisms of Latin squares of order up to $9$.

\section{Constraints related to a Latin square autotopism}

Given a autotopism $\Theta=(\alpha,\beta,\gamma)\in \mathfrak{A}_n$, let $(1)_{\Theta}$ be the set of constraints obtained by adding to $(1)$ the $n^3$ constraints: $$x_{ijk}=x_{\alpha(i)\beta(j)\gamma(k)}, \forall i\in I, j\in J, k\in K. \hspace{1cm} (1.5)_{\Theta}$$

The following results hold:

\begin{theorem} \label{thr1}
There exists a 1-1 correspondence between $LS(\Theta)$ and the set $FS(\Theta)$ of feasible solutions related to a combinatorial optimization problem having $(1)_{\Theta}$ as the set of constraints.
\end{theorem}

{\em Proof.} It is enough to consider the restriction to $LS(\Theta)$ of the correspondence $\Phi$ between $LS(n)$ and $FS(n)$, because then, given $L=(l_{i,j})\in LS(n)$, it is verified that $L\in LS(\Theta)$ if and only if, for all $i,j,k\in [n]$: $l_{i,j}=k \Leftrightarrow l_{\alpha(i),\beta(j)}=\gamma(k)$. But this last condition is equivalent to say that $x_{ijk}=1$ if and only if $x_{\alpha(i)\beta(j)\gamma(k)}=1$. That is to say, $x_{ijk}=x_{\alpha(i)\beta(j)\gamma(k)}$.
\hfill \qed \vspace{0.25cm}

\begin{corollary} \label{crl1} Every feasible solution of $FS(\Theta)$ verifies that $x_{ijk}=0$, for all $i,j,k\in [n]$ such that one of the following assertions is verified:
\begin{enumerate}
\item[a)] $i\in Fix(\alpha), j\in Fix(\beta)$ and $k\not\in Fix(\gamma)$.
\item[b)] $i\in Fix(\alpha), k\in Fix(\gamma)$ and $j\not\in Fix(\beta)$.
\item[c)] $j\in Fix(\beta), k\in Fix(\gamma)$ and $i\not\in Fix(\alpha)$.
\end{enumerate}
\end{corollary}

{\em Proof.} From the conjugacy of rows, columns and symbols in Latin squares, it is enough to consider assertion (a). So, let us consider a feasible solution of $FS(\Theta)$ such that $x_{ijk}=1$, for some $i,j,k\in [n]$ verifying assertion (a). From Theorem \ref{thr1}, there exists an unique $L=(l_{i,j})\in LS(\Theta)$ being equivalent with such a feasible solution. Specifically, it must be $l_{i,j}=k$ and therefore, $k=l_{i,j}=l_{\alpha(i),\beta(j)}=\gamma(l_{i,j})=\gamma(k)$, which is a contradiction, because $k\not\in Fix(\gamma)$. \hfill \qed \vspace{0.25cm}

Let $S_{Fix(\Theta)}$ be the set of triples $(i,j,k)\in [n]^3$ such that one of the assertions of Corollary \ref{crl1} is verified. Since the $\mathbf{l}_{\alpha}^1\cdot \mathbf{l}_{\beta}^1\cdot (n -\mathbf{l}_{\gamma}^1) + \mathbf{l}_{\alpha}^1\cdot \mathbf{l}_{\gamma}^1\cdot (n -\mathbf{l}_{\beta}^1) +
\mathbf{l}_{\beta}^1\cdot \mathbf{l}_{\gamma}^1\cdot (n -\mathbf{l}_{\alpha}^1)$ variables $x_{ijk}$ related to $S_{Fix(\Theta)}$ are all nulls, we can reduce the number of variables of the system $(1)_{\Theta}$ in order to obtain a $1-1$ correspondence between $FS(\Theta)$ and $LS(\Theta)$. Given $s,t\in [n]$, the following sets will be useful:{\small
$$S^{(1,s,t)}_{Fix(\Theta)}=\{i\in [n] \mid (i,s,t)\in S_{Fix(\Theta)}\},\hspace{1cm} S^{(2,s,t)}_{Fix(\Theta)}=\{j\in [n] \mid (s,j,t)\in S_{Fix(\Theta)}\},$$
$$S^{(3,s,t)}_{Fix(\Theta)}=\{k\in [n] \mid (s,t,k)\in S_{Fix(\Theta)}\}.$$}

Moreover, the symmetrical structure given by the autotopism $\Theta$ can also be used to reduce the number of variables of $(1)_{\Theta}$. To see it, let us consider:
$$S_{\Theta}  = \left\{(i,j) \,\mid \, i\in S_{\alpha},j\in \begin{cases} [n],
     \text{ if }i\not \in Fix(\alpha),\\
     S_{\beta}, \text{ if } i\in Fix(\alpha).\end{cases}\right\}$$ as a set of $(\mathbf{n}_{\alpha}-\mathbf{l}_{\alpha}^1)\cdot n + \mathbf{l}_{\alpha}^1\cdot \mathbf{n}_{\beta}$ multi-indices, where $S_{\alpha}=\{c^{\alpha}_{i,1}\mid i\in [\mathbf{n}_{\alpha}]\}$ and $S_{\beta}=\{c^{\beta}_{j,1}\mid j\in [\mathbf{n}_{\beta}]\}$. The following result is verified:

\begin{proposition} \label{prp1} Let $L=(l_{i,j})\in LS(\Theta)$ be such that
all the triples of the Latin subrectangle $R_L=\left\{(i,j,l_{i,j})\,\mid (i,j)\in S_{\Theta}\right\}$ of $L$ are known. Then, all the triples of $L$ are known. Indeed, given $i,j\in [n]$, there exists an unique element $(i_{\Theta},j_{\Theta})\in S_{\Theta}$ such that $l_{i,j}$ can be obtained starting from $l_{i_{\Theta},j_{\Theta}}$.
\end{proposition}

{\em Proof.} Let $(i,j,l_{i,j})\in L$ be such that $i>\mathbf{n}_{\alpha}$ and let $r\in [\mathbf{n}_{\alpha}]$ and $u\in
[\lambda_{r}^{\alpha}]$ be such that $c_{r,u}^{\alpha}=i$. Then, $(\alpha^{1-u}(i),\beta^{1-u}(j))\in S_{\Theta}$, and, therefore, $l_{\alpha^{1-u}(i),\beta^{1-u}(j)}$ is known. Thus,
$l_{i,j}=\gamma^{u-1}(l_{\alpha^{1-u}(i),\beta^{1-u}(j)})$.

Now, let $(i,j,l_{i,j})\in L$ be such that $i\in Fix(\alpha)$ and $j>\mathbf{n}_{\beta}$. Let $s\in [\mathbf{n}_{\beta}]$ and $v\in [\lambda_{s}^{\beta}]$ be such that $c_{s,v}^{\beta}=j$. From the hypothesis, the triple $(i,\ c_{s,1}^{\beta},\ l_{i,
c_{s,1}^{\beta}})$ is known. Thus, $l_{i,j}=\gamma^{v-1}(l_{i,c_{s,1}^{\beta}})$.

The final assertion is therefore an immediate consequence of the election of the cyclic decomposition of $\Theta$. Specifically, it is verified that $(i_{\Theta}, j_{\Theta})=(\alpha^{m_{i,j}}(i),\beta^{m_{i,j}}(j))$, where $m_{i,j}=\min \{t\geq 0 \mid
(\alpha^t(i),\beta^t(j))\in S_{\Theta}\}$.
\hfill \qed \vspace{0.25cm}

Given $i,j,k\in [n]$, let us define $k_{\Theta}=\gamma^{m}(k)$, where $m\in [n]$ is such that $(i_{\Theta}, j_{\Theta})=(\alpha^{m}(i),\beta^{m}(j))\in S_{\Theta}$. Thus, from the cyclic decomposition of $\Theta$, let us observe that $(i_{\Theta},j_{\Theta},k_{\Theta})=(\alpha^t(i)_{\Theta},\beta^t(j)_{\Theta},\gamma^t(k)_{\Theta})$, for all $i,j\in [n]$ and for all $t\in [n]$. The following result holds:

\begin{theorem} \label{thr2} There exists a 1-1 correspondence between $FS(\Theta)$ and the set of feasible solutions $FS'(\Theta)$ of the following system of equations in $d_{\Theta}=((\mathbf{n}_{\alpha}-\mathbf{l}_{\alpha}^1)\cdot n^2 + \mathbf{l}_{\alpha}^1\cdot \mathbf{n}_{\beta}\cdot n)-(\mathbf{l}_{\alpha}^1\cdot \mathbf{l}_{\beta}^1\cdot (n -\mathbf{l}_{\gamma}^1) + \mathbf{l}_{\alpha}^1\cdot \mathbf{l}_{\gamma}^1\cdot (\mathbf{n}_{\beta} -\mathbf{l}_{\beta}^1) +
\mathbf{l}_{\beta}^1\cdot \mathbf{l}_{\gamma}^1\cdot (\mathbf{n}_{\alpha} -\mathbf{l}_{\alpha}^1))$ variables:
$$\begin{cases}\begin{array}{lcr}
         \sum_{i\in [n]\setminus S^{(1,j,k)}_{Fix(\Theta)}}x_{i_{\Theta}j_{\Theta}k_{\Theta}}=1, \forall j,k\in [n]. & \ & (2.1)_{\Theta}\\
         \sum_{j\in [n]\setminus S^{(2,i,k)}_{Fix(\Theta)}}x_{i_{\Theta}j_{\Theta}k_{\Theta}}=1, \forall i,k\in [n]. & \ & (2.2)_{\Theta}\\
         \sum_{k\in [n]\setminus S^{(3,i,j)}_{Fix(\Theta)}}x_{i_{\Theta}j_{\Theta}k_{\Theta}}=1, \forall i,j\in [n]. & \ & (2.3)_{\Theta}\\
         x_{ijk}\in \{0,1\}, \forall (i,j,k)\in S_{\Theta}\times [n] \setminus S_{Fix(\Theta)}. & \ & (2.4)_{\Theta}\end{array}
       \end{cases} \hspace{0.25cm} (2)_{\Theta}
$$
\end{theorem}

{\em Proof.} Let us define the map $\Psi_{\Theta}:FS'(\Theta)\subseteq\mathbb{R}^{d_{\Theta}}\rightarrow FS(\Theta)\subseteq\mathbb{R}^{n^3}$, such that {\small $\Psi_{\Theta}((x_{ijk})_{(i,j,k)\in S_{\Theta}\times [n] \setminus S_{Fix(\Theta)}})=(X_{uvw})_{(u,v,w)\in [n]^3}=\begin{cases} 0, \text{ if } (u,v,w)\in S_{Fix(\Theta)}, \\
  x_{u_{\Theta}v_{\Theta}w_{\Theta}}, \text{ otherwise}.\end{cases}$}. Thus, $\Psi_{\Theta}$ is a 1-1 correspondence between $FS'(\Theta)$ and $FS(\Theta)$. Specifically, from Corollary \ref{crl1} and Proposition \ref{prp1}, equations $(1.1), (1.2)$ and $(1.3)$ and conditions $(1.4)$ in $FS(\Theta)$ are equivalent to $(2.1)_{\Theta},$ $(2.2)_{\Theta}, (2.3)_{\Theta}$ and $(2.4)$ in $FS'(\Theta)$, respectively. Now, let us consider $(x_{ijk})_{(i,j,k)\in S_{\Theta}\times [n] \setminus S_{Fix(\Theta)}}\in FS'(\Theta)$ and $(X_{uvw})_{(u,v,w)\in [n]^3}=\Psi_{\Theta}((x_{ijk})_{(i,j,k)\in S_{\Theta}\times [n] \setminus S_{Fix(\Theta)}})$. Given $u,v,w\in [n]$, it is verified that $X_{uvw}=\begin{cases} 0=X_{\alpha(u)\beta(v)\gamma(w)}, \text{ if } (u,v,w)\in S_{Fix(\Theta)},\\
  x_{u_{\Theta}v_{\Theta}w_{\Theta}}=X_{\alpha(u)\beta(v)\gamma(w)}, \text{ otherwise}.\end{cases}$. Therefore equations $(1.5)_{\Theta}$ are also verified.
\hfill \qed \vspace{0.25cm}

In general, many of the expressions of $(2)_{\Theta}$ are the same equation and so, they are redundant. An immediate consequence of Theorem \ref{thr2} is the following:

\begin{corollary} \label{crl2} $\Psi_{\Theta}^{-1}\, \circ \, \Phi_{|_{LS(\Theta)}}$ is a 1-1 correspondence between $LS(\Theta)$ and $FS'(\Theta)$. \qed \vspace{0.25cm}
\end{corollary}

\section{$0/1$-polytopes related to a Latin square autotopism}

Given a autotopism $\Theta\in \mathfrak{A}_n$, let $A_{\Theta}$ and $A'_{\Theta}$ be the constraint matrices of $(1)_{\Theta}$ and $(2)_{\Theta}$, respectively. Let us define the following $0/1$-polytopes:{\small
$$\mathcal{P}_{LS(\Theta)}=conv\{FS(\Theta)\}=conv\{{\bf x}\in \{0,1\}^{n^3} \mid A_{\Theta}\cdot {\bf x}={\bf e}_{\Theta}\}\subseteq \mathbb{R}^{n^3},$$
$$\mathcal{P}'_{LS(\Theta)}=conv\{FS'(\Theta)\}=conv\{{\bf x}\in \{0,1\}^{n^3} \mid A'_{\Theta} \cdot {\bf x}={\bf e}'_{\Theta}\}\subseteq \mathbb{R}^{d_{\Theta}},$$}
where ${\bf e}_{\Theta}=(1,...,1)^t$ and ${\bf e}'_{\Theta}=(1,...,1)^t$ have $3\cdot n^2 + n^3$ and $3\cdot n^2$ entries, respectively. The following results hold:

\begin{corollary} \label{crl3} Both $0/1$-polytopes, $\mathcal{P}_{LS(\Theta)}$ and $\mathcal{P}'_{LS(\Theta)}$, have $\Delta(\Theta)$ vertices.
\end{corollary}

{\em Proof.}  It is enough to consider the 1-1 correspondences of Theorem \ref{thr1} and Corollary \ref{crl2}.
\hfill \qed \vspace{0.25cm}

\begin{theorem} \label{thr3} $\dim(\mathcal{P}_{LS(\Theta)})=\dim(\mathcal{P}'_{LS(\Theta)})\leq d_{\Theta} - rank(A'_{\Theta})$.
\end{theorem}

{\em Proof.} The inequality is an immediate consequence of the definition of $\mathcal{P}'_{LS(\Theta)}$. Besides, from the definition of $\Psi_{\Theta}$ given in the proof of Theorem \ref{thr2}, it is immediate to see that a set of $m$ affinely vertices of $\mathcal{P}'_{LS(\Theta)}$ induces a set of $m$ affinely vertices of $\mathcal{P}_{LS(\Theta)}$, because we can identify all the coordinates of the first ones in the second ones. So, $\dim(\mathcal{P}'_{LS(\Theta)})\leq\dim(\mathcal{P}_{LS(\Theta)})$.

Now, let $\{V_1,...,V_m\}$ be a set of $m$ affinely independent vertices of $\mathcal{P}_{LS(\Theta)}$, where $V_i=(v_{i,1},...,v_{i,n^3})$, for all $i\in [m]$. From Theorem \ref{thr2}, $V'_i=\Psi_{\Theta}^{-1}(V_i)=(v'_{i,1},...,v'_{i,d_{\Theta}})$ is a vertex of $\mathcal{P}'_{LS(\Theta)}$, for all $i\in [m]$. Let us suppose that there exist $\lambda_1,...,\lambda_m\in \mathbb{R}$, such that $\sum_{i=1}^m \lambda_i=1$ and $\sum_{i=1}^m \lambda_i\cdot V'_i={\bf 0}$. From the definition of $\Psi_{\Theta}$, given $j\in [n^3]$ non corresponding to a triple of $S_{Fix(\Theta)}$, there exists $k\in [d_{\Theta}]$, such that $v_{i,j}=v'_{i,k}$, for all $i\in [m]$. Thus, $\sum_{i=1}^m \lambda_i\cdot V_i={\bf 0}$, which is a contradiction. Therefore, $\dim(\mathcal{P}_{LS(\Theta)})\leq\dim(\mathcal{P}'_{LS(\Theta)})$.
\hfill \qed \vspace{0.25cm}

\begin{theorem} \label{thr4} Let
$(\mathbf{l}_{\alpha},\mathbf{l}_{\beta},\mathbf{l}_{\gamma})$ be
the cycle structure of a Latin square autotopism and let us
consider
$\Theta_1=(\alpha_1,\beta_1,\gamma_1),\Theta_2=(\alpha_2,\beta_2,\gamma_2)\in
\mathfrak{A}_n(\mathbf{l}_{\alpha},\mathbf{l}_{\beta},\mathbf{l}_{\gamma})$.
Then, $\mathcal{P}_{LS(\Theta_1)}$ and $\mathcal{P}_{LS(\Theta_2)}$ are $0/1$-equivalents. Analogously,  $\mathcal{P}'_{LS(\Theta_1)}$ and $\mathcal{P}'_{LS(\Theta_2)}$ are $0/1$-equivalents.
\end{theorem}

{\em Proof.} Let us prove the first assertion, the other case follows analogously. So, since $\Theta_1$ and $\Theta_2$ have the same cycle structure,
we can consider the isotopism $\Theta=(\sigma_1,\sigma_2,\sigma_3)\in \mathcal{I}_n$, where:

\begin{enumerate}
\item[i)] $\sigma_1(c_{i,j}^{\alpha_1})=c_{i,j}^{\alpha_2}$, for all
$i\in [k_{\alpha_1}]$ and $j\in [\lambda_i^{\alpha_1}]$,
\item[ii)] $\sigma_2(c_{i,j}^{\beta_1})=c_{i,j}^{\beta_2}$, for all
$i\in [k_{\beta_1}]$ and $j\in [\lambda_i^{\beta_1}]$,
\item[iii)] $\sigma_3(c_{i,j}^{\gamma_1})=c_{i,j}^{\gamma_2}$, for all
$i\in [k_{\gamma_1}]$ and $j\in [\lambda_i^{\gamma_1}]$.
\end{enumerate}

Let $L\in LS(\Theta_1)$ and $(x_{ijk})_{i,j,k\in [n]}=\Phi(L)$. From \cite{FalconMartinJSC07}, we know that $L\in LS(\Theta_1)$ if and only if $L^{\Theta}\in LS(\Theta_2)$. Thus, if $(X_{ijk})_{i,j,k\in [n]}=\Phi(L^{\Theta})$, then it must be $x_{ijk}=x_{\sigma_1(i)\sigma_2(j)\sigma_3(k)}$, for all $i,j,k\in [n]$. So, the permutation of coordinates $\pi(x_{ijk})=x_{\sigma_1(i)\sigma_2(j)\sigma_3(k)}$ is a 1-1 correspondence between $FS(\Theta_1)$ and $FS(\Theta_2)$, which are the set of vertices of $\mathcal{P}_{LS(\Theta_1)}$ and $\mathcal{P}_{LS(\Theta_2)}$, respectively. Thus, $\pi$ transforms $\mathcal{P}_{LS(\Theta_1)}$ into $\mathcal{P}_{LS(\Theta_2)}$.
\hfill \qed \vspace{0.25cm}

From Theorem \ref{thr4}, the dimension of $\mathcal{P}_{LS(\Theta)}$ and $\mathcal{P}'_{LS(\Theta)}$ only depends on the cycle structure of $\Theta$. Moreover, since rows, columns and symbols have an interchangeable role in Latin squares and since affine independence does not depend on these interchanges, we can suppose that the cycles $\alpha,\beta$ and $\gamma$ of $\Theta$ verify that $\mathbf{n}_{\alpha}\leq \mathbf{n}_{\beta} \leq \mathbf{n}_{\gamma}$. Thus, let us finish this paper by following the classification of all possible cycle structures given in \cite{FalconAC}, in order to show in Tables 1 and 2 the dimensions of all possible polytopes related to any autotopisms of order up to $9$. Specifically, the exact dimension is shown when the set $LS(\Theta)$ is known. As an upper bound we show the difference between $d_{\Theta}$ and $rank(A'_{\Theta})$, which indeed can not be reached, as we can observe in Table 1. As a lower bound, we study the subsets of $LS(\Theta)$ given in \cite{FalconMartinEACA08}.

\begin{table}[ht]
\centering{\tiny
\begin{tabular}{|c|c|c|c|c|c|c|c|c|}\hline
$n$ & $\mathbf{l}_{\alpha}$ & $\mathbf{l}_{\beta}$ &
$\mathbf{l}_{\gamma}$ & $d_{\Theta}$ & $\Delta(\Theta)$ & Lower bound \cite{FalconMartinEACA08} & $\dim(\mathcal{P}'_{LS(\Theta)})$ & $d_{\Theta}$ - rank$(A'_{\Theta})$\\
\hline
2 & (0,1) & (0,1) & (0,1) & 4 & 2 & - & 1 & 1\\
\hline\hline
3 & (0,0,1) & (0,0,1) & (0,0,1) & 9 & 3 & - & 2 & 2\\
\cline{4-4}\cline{6-9}
\ & \ & \ & (3,0,0) & \ & 6 & - & 4 & 4\\
\cline{2-9}
\ & (1,1,0) & (1,1,0) & (1,1,0) & 11 & 4 & - & 3 & 3\\
\hline\hline
4 & \ & \ & (0,2,0,0) & \ & 8 & - & 5 & 7\\
\cline{4-4}\cline{6-9}
\ & (0,0,0,1) & (0,0,0,1) & (2,1,0,0) & 16 & 8 & - & 5 & 8\\
\cline{4-4}\cline{6-9}
\ & \ & \ & (4,0,0,0) & \ & 24 & - & 9 & 9\\
\cline{2-9}
\ & \ & \ & (0,2,0,0) & \ & 32 & - & 12 & 13\\
\cline{4-4}\cline{6-9}
\ & (0,2,0,0) & (0,2,0,0) & (2,1,0,0) & 32 & 32 & - & 13 & 14\\
\cline{4-4}\cline{6-9}
\ & \ & \ & (4,0,0,0) & \ & 96 & - & 15 & 15\\
\cline{2-9}
\ & (1,0,1,0) & (1,0,1,0) & (1,0,1,0) & 19 & 9 & - & 8 & 8\\
\cline{2-9}
\ & (2,1,0,0) & (2,1,0,0) & (2,1,0,0) & 24 & 16 & - & 4 & 7\\
\hline \hline
5 & (0,0,0,0,1) & (0,0,0,0,1) & (0,0,0,0,1) & 25 & 15 & - & 12 & 12\\
\cline{4-4}\cline{6-9}
\ & \ & \ & (5,0,0,0,0) & \ & 120 & - & 16 & 16\\
\cline{2-9}
\ & (1,0,0,1,0) & (1,0,0,1,0) & (1,0,0,1,0) & 29 & 32 & - & 15 & 15\\
\cline{2-9}
\ & (1,2,0,0,0) & (1,2,0,0,0) & (1,2,0,0,0) & 57  & 256 & - & 27 & 28\\
\cline{2-9}
\ & (2,0,1,0,0) & (2,0,1,0,0) & (2,0,1,0,0) & 35 & 144 & - & 15 & 15\\
\hline \hline
6 & \ & \ & (0,0,2,0,0,0) & \ & 72 & - & 19 & 21\\
\cline{4-4}\cline{6-9}
\ & \ & \ & (1,1,1,0,0,0) & \ & 72 & - & 20 & 22\\
\cline{4-4}\cline{6-9}
\ & (0,0,0,0,0,1) & (0,0,0,0,0,1) & (2,2,0,0,0,0) & 36 & 144 & - & 22 & 23\\
\cline{4-4}\cline{6-9}
\ & \ & \ & (3,0,1,0,0,0) & \ & 144 & - & 21 & 23\\
\cline{4-4}\cline{6-9}
\ & \ & \ & (4,1,0,0,0,0) & \ & 288 & - & 24 & 24\\
\cline{4-4}\cline{6-9}
\ & \ & \ & (6,0,0,0,0,0) & \ & 720 & - & 25 & 25\\
\cline{2-9}
\ & (0,0,0,0,0,1) & (0,0,2,0,0,0) & (0,3,0,0,0,0) & 36 & 288 & - & 23 & 23\\
\cline{2-9}
\ & \ & \ & (0,0,2,0,0,0) & \ & 648 & - & 41 & 41\\
\cline{4-4}\cline{6-9}
\ & (0,0,2,0,0,0) & (0,0,2,0,0,0) & (3,0,1,0,0,0) & 72 & 2592 & - & 43 & 43\\
\cline{4-4}\cline{6-9}
\ & \ & \ & (6,0,0,0,0,0) & \ & 25920 & 34 & - & 45\\
\cline{2-9}
\ & (1,0,0,0,1,0) & (1,0,0,0,1,0) & (1,0,0,0,1,0) & 41 & 75 & - & 24 & 24\\
\cline{2-9}
\ & \ & \ & (2,2,0,0,0,0) & \ & 36864 & 37 & - & 63\\
\cline{4-4}\cline{6-9}
\ & (0,3,0,0,0,0) & (0,3,0,0,0,0) & (4,1,0,0,0,0) & 108 & 110592 & 38 & - & 64\\
\cline{4-4}\cline{6-9}
\ & \ & \ & (6,0,0,0,0,0) & \ & 460800 & 27 & - & 65\\
\cline{2-9}
\ & (2,0,0,1,0,0) & (2,0,0,1,0,0) & (2,0,0,1,0,0) & 48 & 768 & - & 25 & 25\\
\cline{2-9}
\ & (2,2,0,0,0,0) & (2,2,0,0,0,0) & (2,2,0,0,0,0) & 88 & 20480 & 20 & - & 44\\
\cline{2-9}
\ & (3,0,1,0,0,0) & (3,0,1,0,0,0) & (3,0,1,0,0,0) & 63 & 2592 & - & 20 & 28\\
\hline \hline
7 & (0,0,0,0,0,0,1) & (0,0,0,0,0,0,1) & (0,0,0,0,0,0,1) & 49 & 133 & - & 30 & 30\\
\cline{4-4}\cline{6-9}
\ & \ & \ & (7,0,0,0,0,0,0) & \ & 5040 & 31 & - & 36\\
\cline{2-9}
\ & (1,0,0,0,0,1,0) & (1,0,0,0,0,1,0) & (1,0,0,0,0,1,0) & 55 & 288 & - & 35 & 35\\
\cline{2-9}
\ & (1,0,2,0,0,0,0) & (1,0,2,0,0,0,0) & (1,0,2,0,0,0,0) & 109 & 42768 & 25 & - & 68\\
\cline{2-9}
\ & (1,1,0,1,0,0,0) & (1,1,0,1,0,0,0) & (1,1,0,1,0,0,0) & 109 & 512 & - & 24 & 52\\
\cline{2-9}
\ & (2,0,0,0,1,0,0) & (2,0,0,0,1,0,0) & (2,0,0,0,1,0,0) & 63 & 4000 & 20 & - & 37\\
\cline{2-9}
\ & (1,3,0,0,0,0,0) & (1,3,0,0,0,0,0) & (1,3,0,0,0,0,0) & 163 & 6045696 & 30 & - & 101\\
\cline{2-9}
\ & (3,0,0,1,0,0,0) & (3,0,0,1,0,0,0) & (3,0,0,1,0,0,0) & 79 & 41472 & 27 & - & 41\\
\cline{2-9}
\ & (3,2,0,0,0,0,0) & (3,2,0,0,0,0,0) & (3,2,0,0,0,0,0) & 131 & 1327104 & 20 & - & 66\\
\hline
\end{tabular}}\caption{Number of vertices and dimensions of polytopes related to $\mathfrak{A}_n$, for $n\leq 7$.}
\end{table}

\begin{table}[ht]
\centering{\tiny
\begin{tabular}{|c|c|c|c|c|c|c|c|}\hline
$n$ & $\mathbf{l}_{\alpha}=\mathbf{l}_{\beta}$ &
$\mathbf{l}_{\gamma}$ & $d_{\Theta}$ & $\Delta(\Theta)$ & Lower bound \cite{FalconMartinEACA08} & $\dim(\mathcal{P}'_{LS(\Theta)})$ & $d_{\Theta}$ - rank$(A'_{\Theta})$\\
\hline
8 & \ & (0,0,0,2,0,0,0,0) & \ & 1152 & - & 43 & 43\\
\cline{3-3}\cline{5-8}
\ & \ & (0,2,0,1,0,0,0,0) & \ & 1408 & - & 44 & 44\\
\cline{3-3}\cline{5-8}
\ & \ & (0,4,0,0,0,0,0,0) & \ & 3456 & 32 & - & 45\\
\cline{3-3}\cline{5-8}
\ & \ & (2,1,0,1,0,0,0,0) & \ & 1408 & - & 44 & 45\\
\cline{3-3}\cline{5-8}
\ & (0,0,0,0,0,0,0,1) & (2,3,0,0,0,0,0,0) & 64 & 3456 & 45 & - & 46\\
\cline{3-3}\cline{5-8}
\ & \ & (4,0,0,1,0,0,0,0) & \ & 3456 & 35 & - & 46\\
\cline{3-3}\cline{5-8}
\ & \ & (4,2,0,0,0,0,0,0) & \ & 8064 & 38 & - & 47\\
\cline{3-3}\cline{5-8}
\ & \ & (6,1,0,0,0,0,0,0) & \ & 17280 & 39 & - & 48\\
\cline{3-3}\cline{5-8}
\ & \ & (8,0,0,0,0,0,0,0) & \ & 40320 & 41 & - & 49\\
\cline{2-8}
\ & \ & (0,0,0,2,0,0,0,0) & \ & 106496 & 38 & - & 85\\
\cline{3-3}\cline{5-8}
\ & \ & (0,2,0,1,0,0,0,0) & \ & 188416 & 43 & - & 86\\
\cline{3-3}\cline{5-8}
\ & \ & (0,4,0,0,0,0,0,0) & \ & 811008 & 36 & - & 87\\
\cline{3-3}\cline{5-8}
\ & \ & (2,1,0,1,0,0,0,0) & \ & 253952 & 34 & - & 87\\
\cline{3-3}\cline{5-8}
\ & (0,0,0,2,0,0,0,0) & (2,3,0,0,0,0,0,0) & 128 & 1007616 & 38 & - & 88\\
\cline{3-3}\cline{5-8}
\ & \ & (4,0,0,1,0,0,0,0) & \ & 712704 & 41 & - & 88\\
\cline{3-3}\cline{5-8}
\ & \ & (4,2,0,0,0,0,0,0) & \ & 2727936 & 35 & - & 89\\
\cline{3-3}\cline{5-8}
\ & \ & (6,1,0,0,0,0,0,0) & \ & 7741440 & 26 & - & 90\\
\cline{3-3}\cline{5-8}
\ & \ & (8,0,0,0,0,0,0,0) & \ & 23224320 & 41 & - & 91\\
\cline{2-8}
\ & (0,1,0,0,0,1,0,0) & (2,0,0,0,0,1,0,0) & 128 & 3456 & 34 & - & 58\\
\cline{3-3}\cline{5-8}
\ & \ & (2,0,2,0,0,0,0,0) & \ & 19008 & 32 & - & 59\\
\cline{2-8}
\ & (1,0,0,0,0,0,1,0) & (1,0,0,0,0,0,1,0) & 71 & 931 & - & 48 & 48\\
\cline{2-8}
\ & \ & (0,2,0,1,0,0,0,0) & \ & 16384 & 17 & - & 112\\
\cline{3-3}\cline{5-8}
\ & (0,2,0,1,0,0,0,0) & (2,1,0,1,0,0,0,0) & 192 & 16384& 18  & - & 113\\
\cline{3-3}\cline{5-8}
\ & \ & (4,0,0,1,0,0,0,0) & \ & 147456 & 19 & - & 114\\
\cline{2-8}
\ & (2,0,0,0,0,1,0,0) & (2,0,0,0,0,1,0,0) & 80 & 19584 & 35 & - & 51\\
\cline{2-8}
\ & \ & (0,4,0,0,0,0,0,0) & \ & - & - & - & 171\\
\cline{3-3}\cline{5-8}
\ & \ & (2,3,0,0,0,0,0,0) & \ & - & - & - & 172\\
\cline{3-3}\cline{5-8}
\ & (0,4,0,0,0,0,0,0) & (4,2,0,0,0,0,0,0) & 256 & - & - & - & 173\\
\cline{3-3}\cline{5-8}
\ & \ & (6,1,0,0,0,0,0,0) & \ & - & - & - & 174\\
\cline{3-3}\cline{5-8}
\ & \ & (8,0,0,0,0,0,0,0) & \ & 828396011520 & 41 & - & 175\\
\cline{2-8}
\ & (2,0,2,0,0,0,0,0) & (2,0,2,0,0,0,0,0) & 152 & 12985920 & 25 & - & 96\\
\cline{2-8}
\ & (2,1,0,1,0,0,0,0) & (2,1,0,1,0,0,0,0) & 152 & 8192 & 15 & - & 74\\
\cline{2-8}
\ & (3,0,0,0,1,0,0,0) & (3,0,0,0,1,0,0,0) & 97 & 388800 & 19 & - & 56\\
\cline{2-8}
\ & (2,3,0,0,0,0,0,0) & (2,3,0,0,0,0,0,0) & 224 & - & - & - & 141\\
\cline{2-8}
\ & (4,0,0,1,0,0,0,0) & (4,0,0,1,0,0,0,0) & 128 & 7962624 & 21 & - & 69\\
\cline{2-8}
\ & (4,2,0,0,0,0,0,0) & (4,2,0,0,0,0,0,0) & 192 & 509607936 & 15 & - & 100\\
\hline \hline
9 & \ & (0,0,0,0,0,0,0,0,1) & \ & 2025 & - & 56 & 56\\
\cline{3-3}\cline{5-8}
\ & \ & (0,0,3,0,0,0,0,0,0) & \ & 7128 & 43 & - & 58\\
\cline{3-3}\cline{5-8}
\ & (0,0,0,0,0,0,0,0,1) & (3,0,2,0,0,0,0,0,0) & 81 & 12960 & 45 & - & 60\\
\cline{3-3}\cline{5-8}
\ & \ & (6,0,1,0,0,0,0,0,0) & \ & 71280 & 47 & - & 62\\
\cline{3-3}\cline{5-8}
\ & \ & (9,0,0,0,0,0,0,0,0) & \ & 362880 & 49 & - & 64\\
\cline{2-8}
\ & \ & (0,0,1,0,0,1,0,0,0) & \ & 15552 & 30 & - & 86\\
\cline{3-3}\cline{5-8}
\ & \ & (0,3,1,0,0,0,0,0,0) & \ & 124416 & 32 & - & 88\\
\cline{3-3}\cline{5-8}
\ & (0,0,1,0,0,1,0,0,0) & (3,0,0,0,0,0,1,0,0) & 162 & 62208 & 33 & - & 88\\
\cline{3-3}\cline{5-8}
\ & \ & (3,3,0,0,0,0,0,0,0) & \ & 1244160 & 24 & - & 90\\
\cline{2-8}
\ & (1,0,0,0,0,0,0,1,0) & (1,0,0,0,0,0,0,1,0) & 89 & 4096& 30  & - & 63\\
\cline{2-8}
\ & \ & (0,0,3,0,0,0,0,0,0) & \ & - & - & - & 170\\
\cline{3-3}\cline{5-8}
\ & (0,0,3,0,0,0,0,0,0) & (3,0,2,0,0,0,0,0,0) & 243 & - & - & - & 172\\
\cline{3-3}\cline{5-8}
\ & \ & (6,0,1,0,0,0,0,0,0) & \ & - & - & - & 174\\
\cline{3-3}\cline{5-8}
\ & \ & (9,0,0,0,0,0,0,0,0) & \ & 948109639680 & 49 & - & 176\\
\cline{2-8}
\ & (1,0,0,2,0,0,0,0,0) & (1,0,0,2,0,0,0,0,0) & 177 & 12189696 & 14 & - & 124\\
\cline{2-8}
\ & (1,1,0,0,0,1,0,0,0) & (1,1,0,0,0,1,0,0,0) & 177 & 69120 & 16 & - & 84\\
\cline{2-8}
\ & (2,0,0,0,0,0,1,0,0) & (2,0,0,0,0,0,1,0,0) & 99 & 438256 & 32 & - & 67\\
\cline{2-8}
\ & (3,0,0,0,0,1,0,0,0) & (3,0,0,0,0,1,0,0,0) & 117 & 3110400 & 13 & - & 73\\
\cline{2-8}
\ & (1,4,0,0,0,0,0,0,0) & (1,4,0,0,0,0,0,0,0) & 353 & - & - & - & 246\\
\cline{2-8}
\ & (3,0,2,0,0,0,0,0,0) & (3,0,2,0,0,0,0,0,0) & 207 & - & - & - & 130\\
\cline{2-8}
\ & (4,0,0,0,1,0,0,0,0) & (4,0,0,0,1,0,0,0,0) & 149 & 199065600 & 18 & - & 87\\
\cline{2-8}
\ & (3,3,0,0,0,0,0,0,0) & (3,3,0,0,0,0,0,0,0) & 297 & - & - & - & 187\\
\hline
\end{tabular}}\caption{Number of vertices and dimensions of polytopes related to $\mathfrak{A}_8$ and $\mathfrak{A}_9$.}
\end{table}

\end{document}